\documentclass[12pt]{article}
\usepackage{graphicx}
\usepackage{amssymb}
\usepackage{amsmath}
\numberwithin{equation}{section}

\begin{document}
\author{Lev Sakhnovich}
\date{December 30, 2006}
\textbf{Rational solutions of KZ equation, case $S_{4}.$}
\begin{center} Lev Sakhnovich  \end{center}
735 Crawford ave., Brooklyn, 11223, New York, USA.\\
 E-mail address: lev.sakhnovich@verizon.net
\begin{center}Abstract \end{center}
We consider the Knizhnik-Zamolodchikov system of linear
differential equations. The coefficients of this system are
generated by elements of the symmetric group $S_{n}$. We
separately investigate the case $S_{4}.$ In this case we solve the
corresponding KZ-equation in the explicit form.
\newpage
\section{Introduction }
1. We consider the Knizhnik-Zamolodchikov differential system (see
[3])
\begin{equation}
\frac{dW}{dz}={\rho}A(z)W,\quad z{\in}C,\end{equation} where
$A(z)$ and $W(z)$ are $n{\times}n$ matrices. We suppose that
$A(z)$ has the form
\begin{equation}
A(z)=\sum_{k=1}^{s}\frac{P_{k}}{z-z_{k}},\quad
s{\geq}2,\end{equation} where $z_{k}{\ne}z_{\ell}$ if
$k{\ne}\ell$. The matrices $P_{k}$ are connected with matrix
representation of the symmetric group $S_{n}$ and are defined by
formulas (2.1)-(2.4). A. Chervov and D. Talalaev formulated the
following interesting conjecture [2].\\ \textbf{Conjecture 1.1.}
\emph{The Knizhnik-Zamolodchikov system $(1.1),(1.2)$ has a
rational fundamental matrix solution when parameter $\rho$ is
integer.}\\ We have proved this conjecture [5] for the case when
$\rho={\pm}1$ In the present paper we separately consider the case
$S_{4}$ and solve the corresponding KZ-equation in the explicit
form.\\
 2.  In a
neighborhood of $z_{k}$ the matrix function $A(z)$
 can be represented in the form
\begin{equation}
A(z)=\frac{a_{-1}}{z-z_{k}}+a_{0}+a_{1}(z-z_{k})+...\end{equation}where
$a_{k}$ are $n{\times}n$ matrices. We investigate the case when
$z_{k}$ is either a regular point of $W(z)$ or a pole. Hence the
following relation
\begin{equation}
W(z)=\sum_{p{\geq}m}b_{p}(z-z_{k})^{p},\quad b_{m}{\ne}0
\end{equation} is true. Here $b_{p}$ are $n{\times}n$ matrices.
We note that $m$ can be negative. \\ \textbf{Proposition
1.1.}(necessary condition, see [4])  \emph{If the solution of
system $(1.1)$ has form $(1.4)$ then $m$ is an eigenvalue of
$a_{-1}$.}\\ We denote by M the greatest integer eigenvalue of the
matrix ${\rho}a_{-1}$. Using relations ((1.3) and (1.4)) we obtain
the assertion.\\ \textbf{Proposition 1.2.}(necessary and
sufficient condition, see [4]) \emph{If the matrix system
\begin{equation}
[(q+1)I_{n}-a_{-1}]b_{q+1}=\sum_{j+\ell=q}{\rho}a_{j}b_{\ell},
\end{equation} where $m{\leq}q+1{\leq}M$ has a solution $b_{m},
b_{m+1},...,b_{M}$ and $b_{m}{\ne}0$ then system $(1.1)$ has a
solution of form $(1.4)$.}\\ We introduce the matrices
\begin{equation}
P_{k}^{-}=I+P_{k},\quad P_{k}^{+}=I-P_{k}.\end{equation}
\section{Explicit solution }
 1. We consider the natural
representation of the symmetric group $S_{n}$ (see [1]). By
$(i;j)$ we denote the permutation which transposes $i$ and $j$ and
preserves all the rest. The
 $n{\times}n$ matrix which corresponds to $(i;j)$ is denoted by
 \begin{equation}
 P(i,j)=[p_{k,\ell}(i,j)], \quad (i{\ne}j).\end{equation}
 The elements $p_{k,\ell}(i,j)$ are equal to zero except
 \begin{equation}p_{k,\ell}(i,j)=1,\quad (k=i,\ell=j);\quad p_{k,\ell}(i,j)=1,\quad
(k=j,\ell=i),\end{equation}
 \begin{equation}p_{k,k}(i,j)=1,\quad
 (k{\ne}i,k{\ne}j).\end{equation} Now we introduce the matrices
 \begin{equation}
 P_{k}=P(1,k+1),
 \quad 1{\leq}k{\leq}n-1.\end{equation}
 and the $1{\times}(n-1)$ vector $e=[1,1, ...,1]$ ,
 and
 $n{\times}n$ matrix \begin{equation}
 S=\left[\begin{array}{cc}
   2-n & e\\
   e^{\tau} & 0 \\
 \end{array}\right].
 \end{equation}
 Using relations (1.38) and (2.1)-(2.4) we deduce that
 \begin{equation}
 T=(n-2)I_{n}+S.\end{equation}
 The eigenvalues of $T$ are defined by the equalities
 \begin{equation} \lambda_{1}=n-1,\quad \lambda_{2}=n-2,\quad
 \lambda_{3}=-1.\end{equation}Hence we have $m_{T}=-1,\quad
 M_{T}=n-1.$ We shall consider the case when $\rho=-1$. In this
 case
the matrix solution $W(z)$ of system (1.1) can be written in the
form
\begin{equation}W(z)=\sum_{k=1}^{s}\frac{L_{k}}{z-z_{k}}+Q(z),\end{equation}
where $L_{k}$ are $n{\times}n$ matrices, $Q(z)$ is a  $n{\times}n$
matrix polynomial of the form
\begin{equation}Q(z)=Q_{-1}z+Q_{0}.\end{equation}
By substituting formulas (2.8) and (2.9) in relation (1.1) we
receive the equalities
\begin{equation}
(I-P_{k})L_{k}=0,\quad 1{\leq}k{\leq}s,\end{equation}
\begin{equation}
\sum_{j{\ne}k}\frac{P_{k}L_{j}+P_{j}L_{k}}{z_{k}-z_{j}}+P_{k}Q(z_{k})=0,
\end{equation}
\begin{equation}Q_{-1}=-TQ_{-1}.\end{equation}
2. Further we consider the  case $S_{4}.$ \\ The vectors $L_{k}$
can be chosen in the following form
\begin{equation}
L_{1}=\mathrm{col}[1,1,-1,-1],\quad
L_{2}=\alpha\mathrm{col}[1,-1,1,-1],\quad
L_{3}=\beta\mathrm{col}[1,-1,-1,1].\end{equation} We note, that
$L_{k}$ which are defined by relation (2.13) satisfy conditions
(2.10). The numbers $\alpha$ and $\beta$ are defined by relations
\begin{equation}\alpha=-\frac{z_{3}-z_{2}}{z_{3}-z_{1}},
\quad \beta=\frac{z_{3}-z_{2}}{z_{2}-z_{1}}.\end{equation} Now we
introduce the matrix polynomial $Q(z)$ (see (2.9)), where
\begin{equation}
Q_{-1}=-[(z_{2}-z_{1})(z_{3}-z_{1})]^{-1}\mathrm{col}[3,-1,-1,-1],\end{equation}
\begin{equation}Q_{0}=[(z_{2}-z_{1})(z_{3}-z_{1})]^{-1}(z_{1}N_{1}+z_{2}N_{2}+z_{3}N_{3}).
\end{equation} Here the vectors $N_{k}$ are defined by relations
\begin{equation}
N_{1}=L_{1}\quad N_{2}=\mathrm{col}[1,-1,1,-1],\quad
N_{3}=\mathrm{col}[1,-1,-1,1].\end{equation} By direct
calculations we see that that matrices $L_{k},\quad (k=1.2,3)$ and
matrix polynomial $Q(z)$ which are defined by relations (2.9) and
(2.13)-(2.17),satisfy the relations (2.10)-(2.12). Hence the
vector function
\begin{equation}
Y_{1}(z)=\sum_{k=1}^{3}\frac{L_{k}}{z-z_{k}}+Q(z) \end{equation}
is the rational solution of system (1.1) when $\rho=-1.$\\
 In order to construct the second rational solution of system (1.1)
 we introduce the vector
 \begin{equation}
 M_{k}=\beta_{k}\mathrm{col}[1,1,1,1],\quad k=1,2,3,\end{equation}
 where the numbers $\beta_{k}$ are defined by the relations
\begin{equation}\beta_{1}=z_{2}-z_{3},\quad \beta_{2}=z_{3}-z_{1},\quad
\beta_{3}=z_{1}-z_{2}.\end{equation} By direct calculations we see
that the numbers $\beta_{k}$ satisfy the equations
\begin{equation}
\frac{\beta_{1}+\beta_{2}}{z_{1}-z_{2}}+\frac{\beta_{1}+\beta_{3}}{z_{1}-z_{3}}=0,
\end{equation}
\begin{equation}
\frac{\beta_{1}+\beta_{2}}{z_{2}-z_{1}}+\frac{\beta_{2}+\beta_{3}}{z_{2}-z_{3}}=0.
\end{equation} It follows from relations (2.19)-(2.22) that
conditions (2.10)-(2.12) will be fulfilled if $L_{k}=M_{k},\quad
Q_{-1}=Q_{0}=0.$ Hence the vector function
\begin{equation}
Y_{2}(z)=\sum_{k=1}^{3}\frac{M_{k}}{z-z_{k}} \end{equation} is the
rational solution of system (1.1) when $\rho=-1.$\\ In order to
construct the next solution $Y_{3}(z)$ of system (1.1) we
introduce the vectors
\begin{equation}m_{1}=\beta_{1}\mathrm{col}[0,0,1,a],\quad
m_{2}=\beta_{2}\mathrm{col}[0,b,0,c],\quad
m_{3}=\beta_{3}\mathrm{col}[0,1,a,0],\end{equation} where
\begin{equation} a=-\frac{\beta_{1}}{\beta_{3}},\quad b=-\frac{\beta_{3}}{\beta_{2}},
\quad c=\frac{\beta_{1}^{2}}{\beta_{2}\beta_{3}}.\end{equation} We
note, that the numbers $\beta_{k}$ as in the previous case are
defined by relations  (2.20). It follows from relations (2.24) and
(2.25) that conditions (2.10)-(2.12) will be fulfilled if
$L_{k}=m_{k},\quad Q_{-1}=Q_{0}=0.$ Hence the vector function
\begin{equation}
Y_{3}(z)=\sum_{k=1}^{3}\frac{m_{k}}{z-z_{k}} \end{equation} is the
rational solution of system (1.1) when $\rho=-1.$\\ In order to
construct the next solution $Y_{4}(z)$ of system (1.1) we
introduce the vectors
\begin{equation}\ell_{1}=\alpha_{1}\mathrm{col}[0,0,1,a],\quad
\ell_{2}=\alpha_{2}\mathrm{col}[0,b,0,c],\quad
\ell_{3}=\alpha_{3}\mathrm{col}[0,d,e,0],\end{equation} where
\begin{equation} a=-\frac{\alpha_{1}}{\alpha_{3}},\quad b=-\frac{\alpha_{3}}{\alpha_{2}}d,
\quad c=\frac{\alpha_{1}^{2}}{\alpha_{2}\alpha_{3}},\end{equation}
\begin{equation}d=-\frac{\alpha_{1}^{2}}{\alpha_{2}\alpha_{3}^{3}}(\alpha_{1}\alpha_{2}
+\alpha_{3}^{2}),\quad
e=1+\frac{\alpha_{1}}{\alpha_{2}}.\end{equation} We note, that the
numbers $\alpha_{k}$ are defined by the relations
\begin{equation}
\alpha_{1}=(z_{3}-z_{2})^{-1},\quad
\alpha_{2}=(z_{1}-z_{3})^{-1},\quad
\alpha_{3}=(z_{2}-z_{1})^{-1}.\end{equation} It follows from
relations (2.27)-(2.30) that conditions (2.10)-(2.12) will be
fulfilled if $L_{k}=\ell_{k},\quad Q_{-1}=Q_{0}=0.$ Hence the
vector function
\begin{equation}
Y_{4}(z)=\sum_{k=1}^{3}\frac{\ell_{k}}{z-z_{k}} \end{equation} is
the rational solution of system (1.1) when $\rho=-1.$\\
\textbf{Proposition 2.1.} \emph{Let us consider the case
$S_{4}$when $\rho=-1.$ In this case the fundamental solution of
system $(1.1)$ has the form
\begin{equation} Y(z)=\sum_{k=1}^{4}c_{k}Y_{k}(z),\end{equation}
where $c_{k}$ are arbitrary constants , the functions $Y_{k}(z)$
are defined by relations $(2.18), (2.23), (2.26)$ and} (2.31).\\
\textbf{Remark 2.1.}The explicit solution for the case $S_{3}$ was
constructed in the paper [5].\\
\begin{center}\textbf{References} \end{center}
1.Burrow M., Representation Theory of Finite Groups, Academic
Press, 1965.\\ 2. Chervov A., Talalaev D., Quantum Spectral
Curves, Quantum Integrable Systems and the Geometric Langlands
Correspondence, arXiv:hep-th/0604128, 2006.\\ 3. Etingof P.I.,
Frenkel I.B., Kirillov A.A. (jr.), Lectures on Representation
Theory and Knizhnik-Zamolodchikov Equations, Amer. Math. Society,
1998.\\ 4. Sakhnovich L.A., Meromorphic Solutions of Linear
Differential Systems, Painleve Type Functions, arXiv:math
CA/0607555,2006\\ 5. Sakhnovich L.A., Rational Solutions of KZ
equation, existence and \\ construction, arXiv.math:-ph/0609067,
2006.\\

\end{document}